\documentclass[11pt]{amsart}
\usepackage{amssymb,amsthm,verbatim, latexsym}
\newtheorem{theorem}{Theorem}[section]

\newtheorem{corollary}[theorem]{Corollary}
\newtheorem{proposition}[theorem]{Proposition}
   \newtheorem{remark}[theorem]{Remark}
\newtheorem{definition}[theorem]{Definition}

\newtheorem{example}[theorem]{Example}

 \newcommand{\nc}{\newcommand}
 \nc{\R}{\ensuremath{\mathbb{R}}}
 \nc{\Z}{\ensuremath{\mathbb{Z}}}
 \nc{\Sim}{\operatorname{S}}
\nc{\D}{\ensuremath{\operatorname{D}}}
\nc{\cl}{\ensuremath{\operatorname{cl}}}
\numberwithin{equation}{section}

\newenvironment{namelist}[1]{%
\begin{list}{}{
     
     \settowidth{\labelwidth}{#1}
      \setlength{\leftmargin}{1.1\labelwidth} }
  }{
\end{list}}
\begin{document}

\title{Dirac's theorem on simplicial matroids}

\author{Raul Cordovil,  Manoel Lemos and Cl\'audia Linhares Sales}

\address{Raul Cordovil\newline
Departamento de Matem\'atica\newline
Instituto Superior T\'ecnico\newline
Av. Rovisco Pais\newline
1049-001 Lisboa, Portugal}
\email{cordovil@math.ist.utl.pt}
\address{Manoel Lemos\newline
Departamento de Matem\'atica\newline
Universidade Federal de Pernambuco\newline
Recife, Pernambuco CEP 50740-540 - Brasil
}
\email{manoel@dmat.ufpe.br}

\address{Cl\'audia Linhares Sales\newline
MDCC,
Departamento de Computa\c c\~ao \newline
Universidade Federal do Cear\'a - UFC\newline
Campus do Pici, Bloco 910\newline
Fortaleza, CE, Brasil}
\email{linhares@lia.ufc.br}

\begin{abstract}
We introduce the notion of $k$-hyper\-clique complexes, i.e.,   the  largest
 simplicial complexes on the set $[n]$  with a fixed $k$-skeleton.
 These simplicial complexes are a  higher-dimen\-sio\-nal analogue of clique (or flag) complexes (case $k=2$) and they
 are a rich new class of simplicial complexes.

  We show that  Dirac's theorem on chordal graphs
has a higher-dimen\-sio\-nal  analogue in which graphs and clique complexes get replaced, respectively,  by  simplicial matroids and $k$-hyper\-clique complexes. We prove also a higher-dimen\-sio\-nal    analogue of  Stanley's reformulation of Dirac's theorem on chordal graphs.
\end{abstract}



\date{\today}

\thanks{First author was partially
supported by FCT (Portugal) through
program POCTI. Second author were partially supported  by CNPq (Grants No. 476224/04-7 and 301178/05-4) and FAPESP/CNPq (Grant No. 2003/09925-5). Third   author  was partially supported  by CNPq (Grants No.
503890/2004-9 and No. 307049/2004-3).}

\maketitle

\section{Introduction and  notations}
 Set $[n] := \{1, 2,\ldots, n \}$.
The
simplicial matroids
$\Sim_{k}^{n}(E)$ on the ground set
$E\subseteq \binom{[n]}{k}$ have been introduced by
 Crapo and Rota
\cite{CraRo, CrapoRota} as one of the six  most important classes  of matroids. These matroids generalize
 graphic matroids:  more precisely
$\Sim_{2}^{n}(E)$ is the  cycle matroid (or graphic
matroid) of the  graph $([n],E)$.  With the aid of Alexander's duality theorem for manifolds applied to simplices, they prove the following beautiful isomorphism
\begin{equation}\label{dual}
\left[\Sim_k^n\left(\binom{[n]}{k}\right)\right]^{*}\simeq \Sim_{n-k}^n\left(\binom{[n]}{n-k}\right), \,X \mapsto [n]\setminus X,
\end{equation}
where $\left[\Sim_k^n\left(\binom{[n]}{k}\right)\right]^{*}$  denotes the dual (or orthogonal) matroid of $\Sim_{n}^k\left(\binom{[n]}{k}\right)$, see \cite[Theorem~11.4]{CrapoRota} or, for an elementary proof (depending only on matrix algebra), \cite[Theorem~6.2.1]{CL}.
In this paper,
we use an equivalent definition of simplicial matroids, see  Definition~\ref{simplicial} below.

We introduce the notion of $k$-hyper\-clique complexes. These simplicial complexes are a natural higher-dimensional analogue of clique (or flag)
complexes (case $k=2$), see Definition~\ref{ksimplicial} below.
The $k$-hyperclique complexes are a rich
  new class of simplicial complexes  of intrinsic interest.
  To  get a better understanding of the structure of $\Sim_k^n(E)$ we attach to it the $k$-hyperclique complex on $[n]$ canonically determined by the family $E$.

 In this paper, we introduce the notion of a strong triangulable simplicial matroid, a higher-dimensional generalization of the notion of a chordal graph. We prove an analogue of Dirac's theorem on chordal graphs
(see Theorem~\ref{Dirac}) using a natural generalization of a perfect sequence of vertices  of a chordal graph (see Theorem~\ref{tri}). We prove also a higher-dimensional  generalization of  Stanley's reformulation of Dirac's theorem on chordal graphs (see Theorem~\ref{main}).

 Let us set some preliminary notation.

\begin{definition}\label{sima}
{\em
An \emph{(abstract) oriented simplicial complex}
 on the set $[n]$ is  a family  $\Delta$  of linear ordered subsets of  $[n] $
 (called the \emph{faces} of $\Delta$)
satisfying the following two conditions. (We identify the  linear ordered
set $\{v_1,v_2,\dotsc, v_m\}$, $v_1<v_2<\dotsm <v_m$, with the symbol $v_1v_2\cdots v_m$.)
\begin{namelist}{xxxx}
\item[$(\ref{sima}.1)$]  Every $v \in [n]$ is a face of $\Delta$;
\item[$(\ref{sima}.2)$] If $F$ is a face of  $\Delta$ and $F' \subset F$,
then $F' $ is also a face of $\Delta$.
\end{namelist}
Given two faces $F'$ and $F=i_1i_2\dotsm i_m$ the ``incidence number" $[F': F]$ is
$$
[F': F]=
\begin{cases}(-1)^{j}&\quad \text{if}\,\,\,  F'=i_1i_2\dotsm{i_{j-1}}{i_{j+1}}\dotsm i_m;\\
0 &\quad \text{otherwhise}.
\end{cases}
$$
  }
   \end{definition}
Let $\mathcal{S}_d(\Delta)$ denote the \emph{$d$-skeleton} of $\Delta$, i.e.,
the  family of
 faces of size $d$
   (or $d$-faces)  of $\Delta$.     A \emph{facet} is a face of $\Delta$, maximal to inclusion.  If nothing in contrary is indicated,  we suppose that
 $\mathcal{S}_1(\Delta)=[n]$.
  \begin{definition}\label{ksimplicial}
  {\em
 Let   $\mathcal{S}_k:=\{F_1,F_2,\dotsc, F_m\}$ be a family of $k$-subsets of $[n]$.
 Let $\langle \mathcal{S}_k \rangle$ be the simplicial complex such that
 $F\subseteq [n]$ is a face of $\langle \mathcal{S}_k \rangle$
 provided that:
 \begin{namelist}{xxxx}
 \item[~$(\ref{ksimplicial}.1)$] $|F|<k$; or
  \item[~$(\ref{ksimplicial}.2)$] if $|F|\geq k$, every $k$-subset of $F$ belongs to
 $\{F_1,F_2,\dotsc, F_m\}$.
\end{namelist}
We say that  $\langle \mathcal{S}_k \rangle$ is the $k$-\emph{hyperclique complex} generated by  the set
 $$\mathcal{S}_k=\{F_1,F_2,\dotsc, F_m\}.$$
 We see $\langle \mathcal{S}_k \rangle$ as an oriented simplicial complex, with the natural orientation induced by $[n]$.
 }
 \end{definition}

Note that $\langle \mathcal{S}_k \rangle$ is the  largest
 simplicial complex $\Delta$ on the set $[n]$ with the fixed $k$-skeleton  $\mathcal{S}_k$.
(In the hypergraph literature, a  family of sets satisfying  Property $(\ref{ksimplicial}.2)$ is said to have the  \emph{Helly dual $k$-property}.)

Throughout this work $ \mathcal{S}_k$ denotes   a family of $k$-subsets of $[n]$ and  let $\langle \mathcal{S}_k \rangle$ denote the corresponding
$k$-hyperclique complex. (So, we have $\mathcal{S}_k=\mathcal{S}_k(\langle \mathcal{S}_k \rangle)$.)
    The paradigm examples of the $k$-hyperclique complexes are the clique (or flag) complexes (the $2$-hyperclique complexes).
  \begin{example}
  {\em Set $\mathcal{S}_k=\{F_1,F_2,\dotsc, F_m\}$.
 If\, $ \bigcap_{i=1}^m F_i\not=\emptyset$  then  $F_1, \dotsc, F_m$ are  facets of $\langle \mathcal{S}_k \rangle$. The other facets of $\langle \mathcal{S}_k \rangle$ are the
 $(k-1)$-subsets of $[n]$ which  are in no  $F_i$.
 For every $k,\, n\geq k\geq 1$, $2^{[n]}$ is a (full) $k$-hyperclique complex.
  }
\end{example}
Let $\mathcal{S}_\ell$ be a subset of $\binom{[n]}{\ell}$, where $n\geq \ell \geq 2$. Let $\langle \mathcal{S}_\ell\rangle$ be the oriented
$\ell$-hyperclique complex determined by $\mathcal{S}_\ell$. Let $\mathbb{F}$ be a field. Consider the two  vector spaces
$\mathbb{F}^{\mathcal{S}_\ell}$ and $\mathbb{F}^{\binom{[n]}{\ell-1}}$ over $\mathbb{F}$.
Let us define the (boundary) map $$\partial_\ell: \mathbb{F}^{\mathcal{S}_\ell} \to \mathbb{F}^{\binom{[n]}{\ell-1}}$$ as the vector space map determined by linearity
 specifying its values in the basis elements:
 \[\partial_\ell F =\sum_{F'\in \binom{[n]}{\ell-1}}[ F':F] F',\]
 for every  $F\in \mathcal{S}_\ell$. By duality
 let us define the (coboundary) map $$\delta^{\ell-1}: \mathbb{F}^{\binom{[n]}{\ell-1}}\to \mathbb{F}^{\mathcal{S}_\ell}$$ as the vector space map determined by linearity
 specifying its values in the basis elements:
\[
\delta^{\ell-1} F' =\sum_{F\in E} [F': F] F\]
for every $F'\in \binom{[n]}{\ell-1}$.
(The symbol $[F': F]$ denotes the incidence number of the faces
$F'$ and $F$ in the oriented simplicial complex   $\langle \mathcal{S}_\ell \rangle$.)
 \begin{definition}\label{simplicial}\cite{CL}
{\em
The \emph{sim\-pli\-cial matroid} $\Sim_{k}^{n}(\mathcal{S}_k)$, on the ground set $\mathcal{S}_k$ and  over the field $\mathbb{F}$, is the  matroid such that
 $$\{X_1, X_2, \dotsc, X_m\}\subseteq \mathcal{S}_k$$ is an
\emph{inde\-pendent  set} iff the
vectors,
$$\partial_k X_1, \partial_k X_2,
\dotsc,
\partial_k  X_m,$$    are linearly
independent in the vector space $\mathbb{F}^{\binom{[n]}{k-1}}$.
}
\end{definition}
\begin{remark}\cite{CLa,CL}
{\em
 Let $(s_{p,q})$ be the matrix whose rows and columns are
 labeled by the sets of $\binom{[n]}{k-1}$ and $\mathcal{S}_k$ respectively, with $s_{p,q}=0$ if $p\not\subseteq q$
 and $s_{p,q}=(-1)^j$ if $q-p={i_j}, q=\{i_1,\dotsc, i_j,\dotsc i_{k}\}$.
The sim\-pli\-cial matroid $\Sim_{k}^{n}(\mathcal{S}_k)$ (over the field $\mathbb{F}$) is the independent matroid of the columns of the
$\{-1,1,0\}$
matrix  $(s_{p,q})$, over the field $\mathbb{F}$. If the matrix $(s_{p,q})$ is not totally unimodular,  the simplicial matroid depends  of the field $\mathbb{F}$. Since the time of Henri Poincar\'e, it is known  that if $ k=2$, the matrix $(s_{p,q})$ is totally unimodular.
The $2$-hyperclique complex $\langle \mathcal{S}_2\rangle$ is the clique  complex of the simple graph
$([n],\mathcal{S}_2)$ and  $\Sim^n_2(\mathcal{S}_2)$ it is   its corresponding  cycle  matroid. So   $\Sim^n_2(\mathcal{S}_2)$ is a regular (or unimodular) matroid, i.e., it is irrespective of the field $\mathbb{F}$.
}
\end{remark}
 If nothing in contrary is said,   the simplicial matroids here considered are over the field $\mathbb{F}$.

For background, motivation, and matroid terminology left undefined here,
see any of the standard references
\cite{CrapoRota, Oxley, Tutte, Welsh} or the  encyclopedic survey    \cite{White1,White2, White3}.  For a description of the
developments on simplicial matroids before 1986, see
\cite{CL}. See also \cite{Bolker}  for an interesting application. For a topological approach    to combinatorics,  see \cite{Bjorner}.  Dirac characterization of chordal graphs (see \cite{D}) is treated extensively
in Chapter 4 of \cite{col}. For an algebraic proof
of Dirac's theorem, see \cite{Hibi}.

\section{Simplicial matroids}\label{app}
The following two propositions are folklore and they are included for completeness. For every vector $v$ of  $\mathbb{F}^E$, where \[v=a_1 e_1+a_2 e_2+\dotsm +a_m e_m\,\, (e_i\in E,\, a_i\in \mathbb{F}^*),\]
 let $\underline{v}:=\{e_1, e_2,\dotsc, e_m\}$ denote  the \emph{support} of  $v$.
\begin{proposition}\label{circuit}
Let $\mathcal{S}_k$ be a subset of $\binom{[n]}{k}$ and
 $\Sim_{k}^{n}(\mathcal{S}_k)$ be the corres\-pon\-ding  simplicial matroid (over the field  $\mathbb{F}$).
 Consider the linear map\newline
 $\partial_k: \mathbb{F}^{\mathcal{S}_k} \to \mathbb{F}^{\binom{[n]}{k-1}}$.
Then
\begin{namelist}{xxxxxxx}
\item[$(\ref{circuit}.1)$] Each circuit of \, $\Sim_{k}^{n}(\mathcal{S}_k)$  has at least $k+1$ elements.
\item[$(\ref{circuit}.2)$] For every $(k+1)$-face   $X$ of $\langle \mathcal{S}_k\rangle
$,
 $\underline{\partial_{k+1} X}$ is a circuit of\, $\Sim_{k}^{n}(\mathcal{S}_k)$.
Each
 circuit with exactly $k+1$ elements is of this type. \qed
\end{namelist}
\end{proposition}
For each $X\in \mathcal{S}_{k+1}(\langle \mathcal{S}_k\rangle)$
we say that   $\underline{\partial_{k+1} X}$ is a \emph{small circuit} of $\Sim_{k}^{n}(\mathcal{S}_k)$.
\begin{proposition}\label{cocircuit}
Let $\mathcal{S}_k$ be a subset of $\binom{[n]}{k}$, $k\geq 2$, and
 $\Sim_{k}^{n}(\mathcal{S}_k)$ be the correspon\-ding  simplicial matroid (over the field  $\mathbb{F}$).
 Consider the linear map\newline $\delta^{k-1}: \mathbb{F}^{\binom{[n]}{k-1} }\to \mathbb{F}^{\mathcal{S}_k}$.
Then
\begin{namelist}{xxxxxxx}
\item[$(\ref{cocircuit}.1)$] The cocircuit space
of \, $\Sim_{k}^{n}(\mathcal{S}_k)$  is generated by the set
 of vectors
$\big\{\delta^{k-1} V\not=0:\, V\in  \binom{[n]}{k-1} \big\}$.
\item[$(\ref{cocircuit}.2)$]
 If non empty, the set\, $\underline{\delta^{k-1} V},\, V\in  \binom{[n]}{k-1}$, is a   union  of cocircuits of $\Sim_{k}^{n}(\mathcal{S}_k)$.  \end{namelist}
\end{proposition}
\begin{proof}
The oriented simplicial complex     $\langle\binom{[n]}{k}\rangle=2^{[n]}$ is  the oriented full\newline
$k$-hyperclique complex.   The matroid $\Sim_k^n(\binom{[n]}{k})$ is the full simplicial matroid
on the ground set $\binom{[n]}{k}$. Consider the linear map
\begin{equation}\label{map}
 \delta^{k-1}: \mathbb{F}^{\binom{[n]}{k-1} }\to \mathbb{F}^{\binom{[n]}{k}}.
  \end{equation}
From   Isomorphism~(\ref{dual}), we know
 that
\[\mathcal{C}^*:=\Big\{\delta^{k-1} V:\, V\in  \binom{[n]}{k-1} \Big\}\]
is a generating set of the  cocircuit space of $\Sim_k^n(\mathcal{S}_k^\prime)$. The linear map  \[\delta^{k-1}: \mathbb{F}^{\binom{[n]}{k-1} }\to \mathbb{F}^{ \mathcal{S}_k}\] is the composition of the map~$(\ref{map})$ and the natural projection \[\iota : \mathbb{F}^{\binom{[n]}{k} }\to \mathbb{F}^{ \mathcal{S}_k}.\]
So,  Assertion $(\ref{cocircuit}.1)$ holds. We know that
$C^*$ is a cocircuit of  $\Sim_k^n(\binom{[n]}{k})$ iff $C^*$
is the support of a non null vector of\, $\mathcal{C}^*$, minimal for inclusion.
Note that $\Sim_k^n(\mathcal{S}_k)^{*}$ is obtained from
$\Sim_k^n(\binom{[n]}{k})^{*}$ by contracting the set $\binom{[n]}{k}\setminus \mathcal{S}_k$.  So, Assertion $(\ref{cocircuit}.2)$ holds.
\end{proof}
Throughout this work  $V, V', V_1, V_2,...$ denote  $(k-1)$-subsets of $[n]$. So, they are $(k-1)$-face of $\langle \mathcal{S}_k\rangle$.
 Let $\langle \mathcal{S}_k\rangle
\!\setminus\!\!\!\setminus V$   denote the $k$-hyper\-clique  complex
   $
   \langle \mathcal{S}_k
\setminus \underline{\delta^{k-1} V}\rangle$, i.e., the $k$-hyper\-clique  complex determined by the set $\mathcal{S}_k
\setminus \underline{\delta^{k-1} V}$.
  Note that, for every pair of      $(k-1)$-faces $V$ and $V'$, we have:
\begin{align*}(\langle \mathcal{S}_k\rangle \!\setminus\!\!\!\setminus V) \!\setminus\!\!\!\setminus V'= \,
 (\langle \mathcal{S}_k\rangle \!\setminus\!\!\!\setminus V')\!\setminus\!\!\!\setminus V=\, \langle \mathcal{S}_k\setminus (\underline{\delta^{k-1} V}\cup  \underline{\delta^{k-1} V'})\rangle.
\end{align*}
\begin{definition}\label{bp}
{\em Let $\Delta_0=\langle \mathcal{S}_k\rangle$ be a $k$-hyper\-clique  complex such that $\Sim_k^n(\mathcal{S}_k)$ has rank
$r$. A sequence $V_{1}, V_{2},\dotsc, V_{r}$ of $(k-1)$-faces of  $\Delta_0$ is said to be \emph{basic linear sequence} when
\begin{equation*} C^*_j:=\underline{\delta^{k-1} V_{j}}\setminus
 \bigcup_{i=1}^{j-1}\underline{\delta^{k-1} V_{i}}
 \end{equation*}
 is a cocircuit of $\Sim_k^n(\mathcal{S}_k(\Delta_{j-1}))$, for $j\in\{1,2,\dotsc, r\}$, where
$$\Delta_{j-1}:=\Delta_{j-2}\setminus\!\!\!\setminus  V_{j-1},\,\, j\in \{2,\dotsc, r\}.$$
 }
\end{definition}
The following result is a corollary of Proposition~\ref{cocircuit}.
\begin{corollary}\label{ell}
Let $\langle\mathcal{S}_k\rangle$ be a $k$-hyperclique complex such that $\Sim_k^n(\mathcal{S}_k)$ has rank
$r$. If\,\,
$\mathcal{V}=(V_1, V_2,\dotsc, V_r)
$
 is a basic linear  sequence of $(k-1)$-faces of
$\langle\mathcal{S}_k\rangle$,
  then
 \begin{equation*}\label{base}
\beta=\{ \delta^{k-1} V_{1}, \delta^{k-1}V_{2},\dotsc, \delta^{k-1}V_{r}\}
\end{equation*}
is a basis of  the cocircuit space  of\, $\Sim_{k}^{n}(\mathcal{S}_k)$.
\end{corollary}
\begin{proof} Suppose  that $\beta$   is a dependent set. Choose a dependent  subset of $\beta$
\begin{equation*} \{\delta^{k-1} V_{i_{1}}, \delta^{k-1} V_{i_{2}}, \dotsc, \delta^{k-1} V_{i_{s}}\},
\end{equation*}
such that $i_1<i_2<\cdots<i_s$ and $s$ is minimum.
Therefore
 \[
\underline{\delta^{k-1} V_{i_{s}}}\subseteq \bigcup_{j=1}^{s-1} \underline{\delta^{k-1} V_{i_{j}}}
\]
and $V_{i_{s}}\not \in \mathcal{V}$,  a contradiction.
As the cocircuit space of  $\Sim_{k}^{n}(\mathcal{S}_k)$ has dimension $r$
the result follows.
\end{proof}
\section{\D-perfect $k$-hyperclique complexes}
In this section we extend to  $k$-hyperclique complexes
the notions of  ``simplicial vertex" and  ``perfect  sequence of vertices", introduced in the Dirac characterization   the clique complexes of chordal graphs, see \cite{D, col}.
\begin{definition}\label{Di}
{\em
Let $
\Delta_0=\langle \mathcal{S}_k\rangle
$ be a
  $k$-hyperclique complex and suppose that the simplicial matroid $\Sim_{k}^{n}(\mathcal{S}_k)$ has rank $r$.
We say that a  $(k-1)$-face $V$ is \emph{simplicial} in $\Delta_0$, if there is exactly one facet  $X$ of
$\Delta_0$ such that  $V\subsetneq X$.
We say that $\Delta_0$  is \emph{\D-perfect} if there is a basic linear sequence of $(k-1)$-faces, $\mathcal{V}=(V_{1}, V_{2},\dotsc, V_{r})$,
such that
every $V_{i}\in \mathcal{V}$  is
simplicial in the $k$-hyperclique complex $\Delta_{i-1}$ where  $$\Delta_{i-1}:=\Delta_{i-2}\!\setminus\!\!\!\setminus V_{i-1},\,\, i\in \{2,\dotsc, r\}.$$
We will call $\mathcal{V}$  a   \emph{\D-perfect sequence
  of $\Delta_0$}.
 }
\end{definition}
 Chordal graphs are an important class of graphs. The following theorem is
one of their fundamental characterizations,  reformulated in our language.
\begin{theorem} {\rm (Dirac's theorem on chordal graphs \cite{D,col})}\label{Dirac}
Let $G=([n], \mathcal{S}_2)$, $\mathcal{S}_2\subseteq \binom{[n]}{2}$ be
a graph and $\langle \mathcal{S}_2\rangle$ be its clique complex.
Then $G$ is chordal if and only if $\langle \mathcal{S}_2\rangle$ is
\D-perfect.\qed
\end{theorem}

\begin{proposition}\label{sim}
Let $V$ be a $(k-1)$-subset of $[n]$.
 If\,  $V$ is  simpli\-cial in the $k$-hyperclique complex $\langle \mathcal{S}_k\rangle$  then $\underline{\delta^{k-1} V}$ is a cocircuit of $\Sim_k^n(\mathcal{S}_k)$.
\end{proposition}
\begin{proof}  From Proposition~\ref{cocircuit} we know that $\underline{\delta^{k-1} V}$
is a   union  of cocircuits of $\Sim_k^n(\mathcal{S}_k)$. Suppose for a contradiction that there are two different  cocircuits $C^*_1$ and $C^*_2$ contained in $\underline{\delta^{k-1} V}$. Choose  elements
$F_1\in C^*_1\setminus C^*_2$ and $F_2\in C^*_2\setminus C^*_1$.  As $V$ is simplicial it follows that $C=\binom{F_1\cup F_2}{k}$ is a circuit of $\Sim_k^n(\mathcal{S}_k)$ and $C\cap C_1^*=\{F_1\}$,
a contradiction to orthogonality. \end{proof}
The reader can easily see
 that the converse of Proposition~\ref{sim}  is not true.
\begin{example}
{\em Set
  \[\mathcal{S}_3=\{ 123,124,125,145,245,136,137,167,367,238,239,289,389\}.\]
Consider the $3$-hyperclique complex $\langle\mathcal{S}_3 \rangle$ on the set $[9]$.
 From Property~(\ref{ksimplicial}.2) we know that
 $\mathcal{S}_4(\langle\mathcal{S}_3 \rangle)= \{1245, 1367, 2389\}$ and
 $\mathcal{S}_5(\langle\mathcal{S}_3 \rangle)=\emptyset$. From Property~(\ref{ksimplicial}.1) we can see that
 the sets of  $2$-faces and $1$-faces   of $\langle\mathcal{S}_3 \rangle$ are respectively
$\mathcal{S}_2(\langle\mathcal{S}_3 \rangle)=\binom{[9]}{2}$ and
$
\mathcal{S}_1(\langle\mathcal{S}_3 \rangle)=\binom{[9]}{1}$.
We can see that the set of facets of $\langle\mathcal{S}_3 \rangle$ is
\begin{align*}
\{& 18,19,26,27,34,35,46,47,48,
 49,56,57,\\
& 58,59,68,69,78,79,123,1245, 1367, 2389\}.
\end{align*}
 Note that $\Sim^9_3(\mathcal{S}_3)$ has rank $10$ and $\langle\mathcal{S}_3 \rangle$ is \D-perfect with the  \D-perfect sequence:
$45,67,89, 15, 14, 16, 17, 28, 29,12$.
}
\end{example}
\begin{proposition}\label{simplicial2}  Let $V$ be a $(k-1)$-subset of $[n]$. Suppose that  $V$  is not a facet
of the $k$-hyperclique complex $\langle \mathcal{S}_k\rangle=\langle F_1, F_2, \dotsc, F_m\rangle$. Then the following two assertions are equivalent:
\begin{namelist}{xxxxxxx}
\item[$(\ref{simplicial2}.1)$]   $V$ is simplicial in\, $\langle \mathcal{S}_k\rangle$;
\item[$(\ref{simplicial2}.2)$] The set $X=\bigcup_{F_i \in \underline{\delta^{k-1} V}} F_i$
is the unique facet of $\langle \mathcal{S}_k\rangle$ containing $V$.
\end{namelist}
\end{proposition}
 \begin{proof}
The implication  $(\ref{simplicial2}.2)\Rightarrow (\ref{simplicial2}.1)$ is  clear.\\
$(\ref{simplicial2}.1)\Rightarrow (\ref{simplicial2}.2)$.
Let $X'$ be the unique facet of $\langle \mathcal{S}_k\rangle$ containing $V$.
Then is clear that $F_i\subseteq X'$ for each  $F_i$ containing $V$.  We conclude that $X\subseteq  X'$ and so $X$ is a face of $\langle \mathcal{S}_k\rangle$.
Suppose, for a contradiction, that $X$ is not a facet of $\langle \mathcal{S}_k\rangle$.
 Then there
 is an $F\in \mathcal{S}_k$ such that
 $F \not \subset X$ but $F\subset X'$. So, for every $x\in F\setminus X$, we know that
 $V\cup x \in\mathcal{S}_k$ and so  $V\cup x \in \underline{\delta^{k-1} V}$. We have the contradiction $F\subset X$. Therefore  $X=X'$.
\end{proof}
\section{Superdense simplicial matroids}
A  matroid $M$ on the ground set $[n]$ and of rank $r$ is called {\em supersolvable}
 if it admits a
maximal chain of modular flats
\begin{equation}\label{sol} \cl(\emptyset)=X_0\subsetneq X_1\subset \dotsm\subsetneq X_{r-1}\subsetneq
X_r=[n].
\end{equation}
 The notion of ``supersolvable lattices''
was introduced and studied by Stanley in  \cite{Stan2}.  For a recent study of supersolvability for chordal binary matroids see \cite{CFS}.
\begin{proposition}\label{solvable}
Let
$\Sim_{k}^{n}(\mathcal{S}_k)$, $k>2$,  be a simplicial matroid.
The  matroid $\Sim_{k}^{n}(\mathcal{S}_k)$ is supersolvable iff it  does not have  circuits.
\end{proposition}
\begin{proof} All the  circuits of $\Sim_{k}^{n}(\mathcal{S}_k)$ have at least $k+1$ elements. So a hyperplane $H$ is modular iff $|\mathcal{S}_k\setminus H|=1$. Indeed if $F, F' \in \mathcal{S}_k\setminus H$, the line $\cl(\{F, F'\})$ cannot  intersect  the hyperplane $H$.  From $(\ref{sol})$ we conclude that  if $\Sim_{k}^{n}(\mathcal{S}_k)$ is supersolvable then it cannot have circuits. The converse is clear.
\end{proof}
So, the notion of supersolvability  is not interesting  for the class of non-graphic  simplicial matroids.
The following definition gives the ``right"  extension of the notion of  supersolvable.
\begin{definition}\label{dense}
{\em
Suppose that $\Sim_{k}^{n}(\mathcal{S}_k)$
has rank $r$.
A hyperplane    $H$ of $\Sim_{k}^{n}(\mathcal{S}_k)$  is said to be \emph{dense} if there is a simplicial   $(k-1)$-face, $V$, of $\langle
\mathcal{S}_k\rangle$ such that: $$H=\mathcal{S}_k\setminus \underline{\delta^{k-1} V}.$$
We say that the simplicial matroid $\Sim_{k}^{n}(\mathcal{S}_k)$ is \emph{superdense}
if it admits a maximal chain of ``relatively dense'' flats
 \begin{equation*}
  \emptyset=X_0 \subsetneq X_1\subsetneq \dotsm\subsetneq X_{r-1}\subsetneq
X_r=\mathcal{S}_k,
\end{equation*}
i.e, such that $X_i$ is a dense hyperplane of\,  $\Sim_{k}^{n}(X_{i+1})$,
 $i\in \{0,1,\dotsc r-1\}$.
}
\end{definition}
  A hyperplane $H$ of $\Sim_{2}^{n}(\mathcal{S}_2)$ is dense
if and only $H$ is modular. Then $\Sim_{2}^{n}(\mathcal{S}_2)$ is superdense
if and only if it is supersolvable.
 So, Theorem~\ref{main} below can be seen as higher-dimensional  generalization of  Stanley's reformulation of Dirac's theorem on chordal graphs, see \cite{Stan2}.
 \begin{theorem}\label{main}
Let $\Delta_0=\langle \mathcal{S}_k\rangle$ be a  $k$-hyperclique
 complex. Then the following two assertions are  equivalent:
\begin{namelist}{xxxxxxxx}
\item[$(\ref{main}.1)$]  $\Delta_0$ is \D-perfect;
\item[$(\ref{main}.2)$]  $\Sim^n_{k}(\mathcal{S}_k)$ is superdense.
\end{namelist}
\end{theorem}
 \begin{proof}
 $(\ref{main}.1)\Rightarrow (\ref{main}.2)$. Let $\mathcal{V}=(V_{1}, V_{2},\dotsc, V_{r})$ be a \D-perfect sequence  of
$\Delta_0$.
 From Proposition~\ref{sim}
we know that the  sets
 \begin{equation*}
 C^*_j:=\underline{\delta^{k-1} V_{j}}\setminus
 \bigcup_{i=1}^{j-1}\underline{\delta^{k-1} V_{i}},\,\,  j\in\{1,2,\dotsc, r\},
 \end{equation*}
 are cocircuits of $\Sim_k^n(\mathcal{S}_k(\Delta_{j-1}))$
 where
$$\Delta_{j-1}:=\Delta_{j-2}\setminus\!\!\!\setminus  V_{j-1},\,\, j\in \{2, 3, \dotsc, r\}.$$ So,  $\mathcal{V}$   determines a
 maximal  chain of  flats of\, $\Sim_{k}^{n}(\mathcal{S}_k)$:
 \begin{equation*}\label{chain}
  \emptyset=X_{0} \subsetneq X_1 \subsetneq  \dotsm\subsetneq  X_{r-1}\subsetneq
X_r=\mathcal{S}_k,
\end{equation*}
where \[X_{r-j}=\mathcal{S}_k(\Delta_{j-1})\setminus C^*_{j},\,\, j=1,\dotsc, r.\]
As  $V_j$ is simplicial in $\Delta_{j-1}$,  we know that $X_{r-j}$ is dense in $\Sim^n_k(\mathcal{S}_k(\Delta_{j-1}))$. So, $\Sim^n_{k}(\mathcal{S}_k)$ is superdense. The proof of the converse part is similar.
 \end{proof}
\section{Triangulable simplicial matroids}
  Now we  introduce a generalization of the notion of  ``triangulable"  for  the classes of
 simplicial matroids. Given a union of  circuits $D$ of $\Sim_k^n(\mathcal{S}_k)$,
 let $\overrightarrow{D}$  denote a vector of $\mathbb{F}^{\mathcal{S}_k}$ whose support is $D$.
 Set $\underline {\overrightarrow{D}}= D$.
 \begin{definition}\label{rigid}
{\em Let $\langle\mathcal{S}_k\rangle=\langle F_1, F_2,\dotsc, F_m\rangle$ be a $k$-hyperclique complex. We say that   $\Sim_{k}^{n}(\mathcal{S}_{k})$ (over the field $\mathbb{F}$)  is \emph{triangulable} provided that the  vector family
\[\{\partial_{k+1} X:\, X \in \mathcal{S}_{k+1}(\langle\mathcal{S}_k\rangle)\}\]
spans the  circuit space.\\
 Moreover, when
 generators $\partial_{k+1}X_1, \partial_{k+1}X_2,\dots, \partial_{k+1} X_{m'}$ can be chosen such
that, for every circuit $C$, there are non-null scalars $a_j\in \mathbb{F}^*$ such that  \[\overrightarrow{C}=\sum_{j=1}^s a_j  \partial_{k+1} X_{i_j}\,\,\, \,\,\,\text{and}\,\,\,\,\,\,
\bigcup_{F_\ell \in \underline{C}} F_\ell = \bigcup_{i=1}^{s} X_{i_j} \,\,\, \,\,\,\text{where}\,\,\,\,\,\, X_{i_j}\in\{ X_1,\dotsc, X_{m'}\}\]
we say that
$\Sim_{k}^{n}(\mathcal{S}_{k})$  is \emph{strongly triangulable}.
}
\end{definition}
Note that we can replace in Definition~\ref{rigid}  the circuit $C$ by a union of circuits $D$.
  It is clear that a simple graph $([n], \mathcal{S}_2)$ is chordal iff $\Sim_2^n(\mathcal{S}_2)$ is strongly triangulable.  Theorem~\ref{tri} is the possible generalization of Dirac's theorem on chordal graphs (see Theorem~\ref{Dirac} above). Indeed, if $k>2$, the converse of Theorem~\ref{tri}   is not true,
  see the remarks following the theorem.
  \begin{theorem}\label{tri}
  Let $\Delta_0=\langle \mathcal{S}_k\rangle=\langle F_1, F_2,\dotsc, F_m\rangle$ be a $k$-hyperclique complex. If   $\Delta_0$ is  \D-perfect,
  then $\Sim_{k}^{n}(\mathcal{S}_k)$ is strongly triangulable.
   \end{theorem}
  \begin{proof}
    The proof is algorithmitic.  Let $\mathcal{V}=(V_1,\dotsc, V_r)$ be a \D-perfect sequence.
   Let $D$ be a union of circuits of $\Sim_{k}^{n}(\mathcal{S}_k)$.
Let $V_i$ the first $(k-1)$-face of $\mathcal{V}$ contained in an
element of $D$. From the definitions
 we know that $V_i$ is a simplicial $(k-1)$-face of $\Delta_{i-1}$
 and $D$  is a union of circuits of $\Sim_{k}^{n}(\mathcal{S}_k(\Delta_{i-1}))$, where
 \[
 \Delta_{i-1}=\Delta_{i-2}\setminus\!\!\!\setminus V_{i-1},\,\, i\in \{2,\dotsc, r\}.\]
From Proposition~\ref{sim}  we  know  that
  \begin{equation*} C^*_j:=\underline{\delta^{k-1} V_{j}}\setminus
 \bigcup_{i=1}^{j-1}\underline{\delta^{k-1} V_{i}}
 \end{equation*}
 is a cocircuit of $\Sim_{k}^{n}(\mathcal{S}_k(\Delta_{i-1}))$.
Set $D\cap C^*_j=\{F_{i_1}, F_{i_2},\dotsc, F_{i_{h}}\}$
 and consider the family of  vectors of $\mathbb{F}^{\mathcal{S}_k}$
\[ \Big\{\overrightarrow{C_{s}}= \partial_{k+1} (F_{i_1}\cup F_{i_s}),\,\,\, s=2, \dotsc, h\Big\}.\]
Express a vector $\overrightarrow{D}$ of support $D$ in the canonical basis, say \[\overrightarrow{D}= a_{i_1}F_{i_1}+ a_{i_2}F_{i_1}+\dotsm + a_{i_h}F_{i_h}+
a_{i_{h+1}}F_j+\dotsm +a_{i_m} F_{i_m},\]
where $a_{i_\ell }\in \mathbb{F}^*, \ell=1,\dotsc,h$, $a_{i_\ell }\in \mathbb{F}, \ell=h+1,\dotsc, m$ and
$\{F_{i_1},\dotsc, F_{i_m}\}=\mathcal{S}_k$.
For every $s\in\{2,3,\dots,h\}$,  it is possible to choose $b_s\in \mathbb{F}^*$ such that $F_{i_s}$ does not belong to
the support of $b_s\overrightarrow{C_s}+\overrightarrow{D}$. As $(C_s\cap D)\cap
C_j^*=\{F_{i_1},F_{i_s}\}$ it follows
that $F_{i_2},F_{i_3},\dots,F_{i_h}$ does not belong to the support of
$$\overrightarrow{D'}:=\overrightarrow{D}+b_2\overrightarrow{C_2}+b_3\overrightarrow {C_3}+\cdots+b_h\overrightarrow{C_h}.$$
The dependent set  $D'$
   is a   union of circuits and  $D'\cap
C_j^*\subseteq\{F_{i_1}\}$. So, by orthogonality we have $D'\cap C_j^*=\emptyset$.
 Note that
\begin{namelist}{xxxx}
\item[$(i)$] For every $V_j\in \mathcal{V}$, $1\leq j \leq i$, no element of $D'$ contain  $V_j$;
 \item[$(ii)$]  $$\bigcup_{F_\ell \in D} F_{\ell}=
\bigcup_{F_{\ell'} \in  \bigcup_{s=2}^h C_s\cup D'} F_{\ell'}.$$
   \end{namelist}
    Replace $D$ by the   set $D'$ and apply the same arguments. From
   $(i)$ we know that  the algorithm   finish.
 It finishes only if  $D'$ is a  small circuit. So the theorem follows.
 \end{proof}
If $k>2$,   the converse of Theorem~\ref{tri}   is not true. Indeed
 consider the triangulation  of a projective plane
 \[F_1=124,\, F_2=126,\, F_3=134,\, F_4=135,\, F_5=165, \]
 \[F_6=235,\, F_7=236,\, F_8=245,\, F_9=346,\, F_{10}=456.\]
 Consider the
  $3$-hyperclique complex
$\langle \mathcal{S}_k\rangle=\langle  F_1, F_2,\dotsc, F_{10}\rangle$ on the set $[6]$.
The simplicial matroid
$\Sim_3^n(\mathcal{S}_k)$, over a field $\mathbb{F}$ of characteristic different of $2$, does not have circuits and then it is (trivially) strongly triangulable. Every $2$-face of a $F_i$  is contained in
 another $F_j$, $j\in \{1,\dotsc, 10\}, j\not= i$.  The facets of $\langle \mathcal{S}_k\rangle$
 are the sets $F_1, F_2,\dotsc, F_{10}$ and all the 2-faces of $\langle \mathcal{S}_k\rangle$
 not contained in an $F_i$. Then $\langle \mathcal{S}_k\rangle$ does not contain simplicial $2$-faces and it is not \D-perfect.

  We remark that the  cycle matroid of a non chordal graph can
be  triangulable. More generally we have:
\begin{proposition}\label{merda}
For any $n, k,\, n-3\geq k\geq 2$, there is a  $k$-hyperclique complex
$\langle \mathcal{S}_k \rangle$ such that:
\begin{namelist}{xxxxxxx}
\item[$(\ref{merda}.1)$] The simplicial matroid $\Sim_k^n(\mathcal{S}_k)$ is triangulable but not strongly triangulable;
\item[$(\ref{merda}.2)$] $\langle \mathcal{S}_k\rangle$ does not contain a simplicial $(k-1)$-face.
\end{namelist}
 \end{proposition}
\begin{proof}
 Let $\langle \mathcal{S}_k\rangle$ be the $k$-hyperclique complex where
\begin{align*}\mathcal{S}_k= &\binom{12\dotsm (k+1)}{k}\bigcup \binom{23\dotsm (k+2)}{k}\setminus 23\dotsm (k+1)\bigcup\\
&\bigcup _{i=1}^{k+1}\binom{12\dotsm\widehat{i}\dotsm (k+1)\,n}{k}
\bigcup\\
& \bigcup _{j=2}^{k+2}\binom{23\dotsm\widehat{j}\dotsm (k+2)\,n}{k}.
\end{align*}
The simplicial matroid $\Sim_{k}^{n}(\mathcal{S}_k)$
has $2k$ small circuits,
\[C_i:=\binom{12\dotsm\widehat{i}\dotsm (k+1)\,n}{k},\,\,\, i\in \{1,2,\dotsc,k+1\},
\]
\[C_j:=\binom{23\dotsm\widehat{j}\dotsm  (k+2)\,n}{k},\,\,\, j\in \{2,3, \dotsc,k+2\}.
\]
The set \[C:=\binom{12\dotsm (k+1)}{k}\bigcup \binom{23\dotsm (k+2)}{k}\setminus 23\dotsm (k+1)\]
is a circuit, symmetric difference of all the $2k$ small circuits.
So,  the simplicial matroid $\Sim_k^n(\mathcal{S}_k)$ over a field of characteristic 2   is triangulable  but not strongly triangulable. The reader can check that do not exist simplicial $(k-1)$-faces  in $\langle \mathcal{S}_k\rangle$.
\end{proof}


\begin{thebibliography}{99}
\bibitem{Bjorner} A.  Bj\"orner:~ Topological methods, in: R. Graham, M. Gr\"otschel, L. Lov\'asz
 (Eds.), Handbook of Combinatorics, \emph{North-Holland},
 Amsterdam, 1995, pp. 1819--1872.


\bibitem{Bolker}
E.\ D. Bolker:~
Simplicial geometry and transportation polytopes.
\emph{Trans. Amer. Math. Soc.}, \textbf{217} (1976), 121--142.

\bibitem{CLa}  Cordovil, R. and Las Vergnas, M.:~ ``G\'eometries simpliciales unimodulaires".
(French) {\em Discrete Math.} {\bf 26} (1979), no. 3, 213\---217.


  \bibitem{CL} R. Cordovil  and B. Lindstr\"om:~ Simplicial matroids, in \cite{White2}.




\bibitem{CFS}R. Cordovil, D. Forge and S. Klein:~ How is a chordal graph like a supersolvable binary matroid? \emph{Discrete Mathematics}, \textbf{288} (2004), 167--172.

\bibitem{CraRo} H.\  H. Crapo and G.-C. Rota:~
Simplicial geometries. Combinatorics.
(Proc. Sympos. Pure Math., Vol. XIX, Univ. California, Los Angeles,
Calif., 1968),
pp. 71--75. \emph{Amer. Math. Soc.,} Providence, R.I., 1971.


\bibitem{CrapoRota}
H.\ H. Crapo and G.-C. Rota:~
On the foundations of combinatorial theory: Combinatorial geometries.
Preliminary edition. \emph{The M.I.T. Press,} Cambridge, Mass.-London, 1970.



\bibitem{D} G.\ A. Dirac:~ On rigid circuit graphs. \emph{Abh.\ Math.\ Sem.\ Univ.\ Hamburg}, \textbf{38} (1961), 71--76.




\bibitem{col} M.\ C. Golumbic:~  Algorithmic graph theory
and perfect graphs. Second Edition. Annals of Discrete Mathematics, \textbf{57}. \emph{Elsevier/North--Holland Press}, New York,  2004.

\bibitem{Hibi} J. Herzog, T. Hibi and  X. Zheng:~ Dirac's theorem on chordal graphs and Alexander duality.
\emph{Europ. J. Combinatorics}, \textbf{25} (2004), 949--960.

\bibitem{Oxley}
J.\ G. Oxley:~  Matroid theory.
Oxford Science Publications.
\emph{The Clarendon Press, Oxford University Press,} New York, 1992.

\bibitem{Stan2} R. \ P. Stanley :~  Supersolvable lattices.
\emph{Algebra Universalis} \textbf{2} (1972), 197\---217.

\bibitem{Tutte} W. \ T. Tutte:~ Introduction to the Theory of Matroids, \emph{American Elsevier,} New York, 1971.

\bibitem{Welsh}
D.\ J.\ A. Welsh:~
Matroid theory.
London Math. Soc. Monographs \textbf {8}
\emph{Academic Press,} London-New York, 1976.

\bibitem{White1}
N. White  (Ed.):~
Theory of matroids.
Encyclopedia Math. Appl., \textbf {26}.
\emph{Cambridge University Press,} Cambridge, 1986.

\bibitem{White2}
N. White  (Ed.):~
Combinatorial geometries.
Encyclopedia Math. Appl., \textbf {29}.
\emph{Cambridge University Press,} Cambridge, 1987.

\bibitem{White3}
N. White  (Ed.):~
Matroid applications.
Encyclopedia Math. Appl.,  \textbf {40}.
\emph{Cambridge University Press,} Cambridge, 1992.


\end{thebibliography}
\end{document}